\documentclass[11pt]{amsart}

\usepackage{geometry}
\usepackage{todonotes}

\usepackage{amsfonts, amssymb, amscd}
\numberwithin{equation}{section}

\usepackage[symbol]{footmisc}

\usepackage{bm}
\usepackage{verbatim}
\usepackage{mathrsfs}
\usepackage{graphicx}
\usepackage{tikz-cd}
\usepackage{subcaption}
\usepackage{listings}
\usepackage{subfiles}
\usepackage[toc,page]{appendix}
\usepackage{mathtools}
\usepackage{comment}
\usepackage{enumerate}
\usepackage{enumitem}
\usepackage[all]{xy}

\usepackage{graphicx}
\graphicspath{{images/}}

\usepackage{appendix}
\usepackage{hyperref}
\hypersetup{
    colorlinks=true,
    citecolor=red,
    linkcolor=blue,
    filecolor=magenta,      
    urlcolor=red,
}
\newcommand{\bQ}{\mathbb{Q}}

\newcommand{\Rr}{\mathbb{R}}

\newcommand{\vol}{\operatorname{vol}}

\newcommand{\tmld}{{\operatorname{tmld}}}

\newcommand{\lcm}{\operatorname{lcm}}

\newcommand{\lct}{\operatorname{lct}}

\newcommand{\Supp}{\operatorname{Supp}}

\newcommand{\mult}{\operatorname{mult}}

\newcommand{\Ii}{\Gamma}

\newtheorem{thm}{Theorem}[section]

\newtheorem{cor}[thm]{Corollary}
\newtheorem{lem}[thm]{Lemma}

\newtheorem{claim}[thm]{Claim}

\theoremstyle{definition}
\newtheorem{defn}[thm]{Definition}

\theoremstyle{definition}
\newtheorem{rem}[thm]{Remark}

\theoremstyle{definition}

\begin{document}

\title{Remark on complements on surfaces}

\author{Jihao Liu}

\address{Department of Mathematics, Northwestern University, 2033 Sheridan Rd, Evanston, IL 60208}
\email{jliu@northwestern.edu}

\date{\today}

\begin{abstract}
We give an explicit characterization on the singularities of exceptional pairs in any dimension. In particular, we show that any exceptional Fano surface is $\frac{1}{42}$-lc. As corollaries, we show that any $\mathbb R$-complementary surface $X$ has an $n$-complement for some integer $n\leq 192\cdot 84^{128\cdot 42^5}\approx 10^{10^{10.5}}$, and Tian's alpha invariant for any surface is $\leq 3\sqrt{2}\cdot 84^{64\cdot 42^5}\approx 10^{10^{10.2}}$. Although the latter two values are expected to be far from being optimal, they are the first explicit upper bounds of these two algebraic invariants for surfaces. 
\end{abstract}

\maketitle

\section{Introduction}

We work over the field of complex number $\mathbb C$. 

Birkar famously proved the boundedness of $n$-complements for $\mathbb R$-complementary varieties and pairs with hyperstandard coefficients \cite{Bir19}, which was later generalized to arbitrary DCC coefficients \cite{HLS19} and arbitrary coefficients \cite{Sho20} under milder conditions. It is interesting to ask whether we can give an explicit bound of $n$, as such explicit bound is expected to be useful for the moduli of log surfaces (cf. \cite{AL19,Kol94,Kol13}) and threefold minimal log discrepancies \cite{HLL22}. \cite{Sho00} shows that an $\mathbb R$-complementary surface pair $(X,B)$ is $n$-complementary for some $n\in\{1,2,3,4,6\}$ when $B$ has standard coefficients and $(X,B)$ is not exceptional, but the question remained open in general for surfaces. In this paper, we provide an explicit upper bound of $n$ for surfaces.

\begin{thm}\label{thm: surface precise complement}
Let $X/Z\ni z$ be an $\mathbb R$-complementary surface. Then $X/Z\ni z$ has an $n$-complement for some $n\leq 192\cdot 84^{128\cdot 42^5}$. In particular, if $Z$ is a point, then $h^0(-nK_X)>0$.
\end{thm}

The key ingredient of the proof of Theorem \ref{thm: surface precise complement} is the following result which provides an explicit characterization of the singularities of exceptional pairs in any dimension. Recall that $\lct(d,\Ii)$ is the set of lc thresholds for effective Weil divisors with respect to pairs of dimension $d$ with coefficients in $\Ii$.

\begin{thm}\label{thm: singularity of exceptional fano varieties}
Let $d$ be a positive integer and $\Ii\subset [0,1]$ a DCC set.  Let
\begin{align*}
   \epsilon_1(d,\Ii):=\inf\Bigg\{1-t\Biggm|
    \begin{array}{r@{}l}t<1,
    \text{there exists a pair }(W,\Delta+t\Psi)\text{ of dimension } d\text { such that }\\ (W,\Delta+t\Psi)\text{ is lc}, K_W+\Delta+t\Psi\equiv 0,
     \Delta\in\Ii,\text{and }0\not=\Psi\in\mathbb N^+
    \end{array}\Bigg\},
    \end{align*}
$$\epsilon_2(d,\Ii):=\inf\{1-t\mid t<1, t\in\lct(d,\Ii)\},$$
and $\epsilon(d,\Ii):=\min\{\epsilon_1(d,\Ii),\epsilon_2(d,\Ii)\}$. Then for any exceptional Fano type pair $(X,B)$ of dimension $d$ such that $B\in\Ii$ and any $0\leq G\sim_{\mathbb R}-(K_X+B)$, $(X,B+G)$ is $\epsilon(d,\Ii)$-lc. In particular, $(X,B)$ is $\epsilon(d,\Ii)$-lc.
\end{thm}

We have the following corollary, which implies Theorem \ref{thm: surface precise complement}.

\begin{cor}\label{cor: 1/42 exceptional fano type surface}
Exceptional Fano type surfaces are $\frac{1}{42}$-lc.
\end{cor}

With Corollary \ref{cor: 1/42 exceptional fano type surface}, we also provide an explicit upper bound of Tian's $\alpha$-invariant for surfaces:

\begin{cor}\label{cor: upper bound alpha invariant surface}
Tian's $\alpha$-invariant for any surface is $\leq 3\sqrt{2}\cdot 84^{64\cdot 42^5}$ (when it is well-defined).
\end{cor}

Although the bounds in Theorem \ref{thm: surface precise complement} and Corollary \ref{cor: upper bound alpha invariant surface} are expected to be far from being optimal, these are the first precise upper bounds of these two algebraic invariants for surfaces. Similar topics and alternative directions include the estimation of the lower bound $n$ (cf. \cite{ETW21,TW21}), the boundedness of the the anti-canonical volume of Fano varieties (cf. \cite{Pro05, Pro07, CC08, Pro10, Che11, Pro13, CJ16, CJ20, JZ21, Jia21b, JZ22, Bir22}), estimation of $(\epsilon,n)$-complement \cite{CH21}, the explicit M\textsuperscript{c}Kernan-Shokurov conjecture \cite{HJL22}, precise bounds of mlds \cite{Jia21a,LX21,LL22}, etc.

\medskip

\noindent\textbf{Acknowledgement}. The author thanks Guodu Chen, Jingjun Han, Kai Huang,  Yuchen Liu, Yujie Luo, Fanjun Meng, Lingyao Xie, Qingyuan Xue, and Junyan Zhao for useful discussions. The author thanks Jingjun Han, Chen Jiang, and Yuchen Liu for proposing the question to the author in different occasions, and thanks Yuchen Liu for the reference \cite{BJ20}. The author thanks Jingjun Han, Chen Jiang, and Dae-Won Lee for comments on the first version of this paper. He thanks the referee for useful suggestions and revising his paper.

\medskip

\noindent\textbf{Postscript}. After the first version of this paper: 1) Totaro \cite{Tot22} conjectured that the smallest Tian's alpha invariant for del Pezzo surfaces is equal to $\frac{21}{2}$, given by $X_{154}\subset\mathbb P(77,45,19,14)$ 2) The author and Shokurov \cite{LS23} prove that $\epsilon(2,\{0\})=\frac{1}{13}$. This result allows us to get a better explicit bounds of the $n$-complements and the $\alpha$-invariants. Nevertheless, in order to make the paper self-contained, we will not use any results in \cite{LS23}.

\section{Preliminaries}

We adopt the standard notation and definitions in \cite{KM98,BCHM10} and will freely use them. For the notation of (relative) pairs $(X/Z\ni z,B)$ and complements, we refer the reader to \cite{CH21}.

\begin{defn}
Let $(X/Z\ni z,B)$ be an $\Rr$-complementary pair. We say that $(X/Z\ni z,B)$ is \emph{exceptional} if $(X/Z\ni z,B+G)$ is klt for any $G\geq 0$ such that $K_X+B+G\sim_{\mathbb R}0$ over a neighborhood of $z$.
\end{defn}

\begin{defn}
Let $d$ be a positive integer and $\Ii\subset [0,1]$ a set. We let 
\begin{align*}
   \epsilon_1(d,\Ii):=\inf\Bigg\{1-t\Biggm|\begin{array}{r@{}l} t<1,
    \text{there exists a pair }(X,B+tC)\text{ of dimension } d\text { such that }\\ (X,B+tC)\text{ is lc}, K_X+B+tC\equiv 0,
     B\in\Ii,\text{and }0\not=C\in\mathbb N^+
    \end{array}\Bigg\},
    \end{align*}
$$\epsilon_2(d,\Ii):=\inf\{1-t\mid t<1, t\in\lct(d,\Ii)\},$$
and $\epsilon(d,\Ii):=\min\{\epsilon_1(d,\Ii),\epsilon_2(d,\Ii)\}$.  By \cite[Theorem 1.5]{HMX14}, $\epsilon(d,\Ii)>0$ when $\Ii$ is DCC.
\end{defn}

\begin{rem}
Usually, $\epsilon_1(d,\Ii)<\epsilon_2(d,\Ii)$, (e.g. when $\Ii=D(\Ii)$ \cite[Lemma 11.2, Proposition 11.5]{HMX14}) and in such cases, $\epsilon(d,\Ii)=\epsilon_1(d,\Ii)$.
\end{rem}

\begin{lem}[{\cite[5.3 Theorem, 5.4 Theorem]{Kol94}}]\label{lem: surface lct 1-gap}
Let $\Ii:=\{1-\frac{1}{n}\mid n\in\mathbb N^+\}\cup\{1\}$. Then:
\begin{enumerate}
    \item $\epsilon_2(2,\Ii)=\epsilon_2(2,\{0\})=\frac{1}{6}$.
    \item $\epsilon_1(2,\Ii)=\frac{1}{42}\leq\epsilon_1(2,\{0\})$.
\end{enumerate}
\end{lem}

\begin{rem}\label{rem: gloabl local correspondenc mld gap}
It was expected that $\epsilon_1(2,\{0\})=\frac{1}{13}$ (cf. \cite[Notation 4.1]{AL19}, \cite[40]{Kol13}). After the first version of this paper, the author and Shokurov prove this result in \cite{LS23}. It is interesting to ask whether $\epsilon_1(d,\Ii)$ is equal to the $1$-gap of mlds for pairs with coefficients in $\Ii$ in dimension $d+1$ (cf. \cite{Jia21a,LX21,LL22}).
\end{rem}

\begin{rem}
When $\Ii=\{1-\frac{1}{n}\mid n\in\mathbb N^+\}\cup\{1\}$ is the standard set, $\epsilon(d,\Ii)\leq\frac{1}{N_{d+2}-1}$ by considering the example $(\mathbb P^d,\sum_{i=1}^{d+1}(1-\frac{1}{N_i})H_i+(1-\frac{1}{N_{d+2}-1})H_{d+2}))$, where $\{N_i\}_{i=1}^{+\infty}$ is the Sylvester sequence $2,3,7,43,\dots$ and $H_i$ are general hyperplanes of degree $1$. It is also expected that $\epsilon(d,\Ii)=\frac{1}{N_{d+2}-1}$ \cite[8.16]{Kol97}. 
\end{rem}

\begin{lem}[cf. {\cite[Proof of Lemma 3.7]{AM04},~\cite[After Theorem A]{Lai16}, \cite[Lemma 2.2]{Bir22}}]\label{lem: index epsilonlc fano surface}
Let $\epsilon$ be a positive real number and $X$ an $\epsilon$-lc Fano surface. Then $IK_X$ is Cartier for some positive integer
$I\leq 2\left(\frac{2}{\epsilon}\right)^{\frac{128}{\epsilon^5}}.$
\end{lem}

We will frequently use the following result to run minimal model programs:

\begin{thm}[{\cite[Corollary 1.3.2]{BCHM10}}]\label{thm: fano type mds}
Fano type varieties are Mori dream spaces. In particular, for any Fano type variety $X$ and any $\Rr$-Cartier $\Rr$-divisor $D$ on $X$, any sequence of $D$-MMP terminates with either a good minimal model or a Mori fiber space.
\end{thm}

\section{The non-exceptional case}

\begin{lem}\label{lem: complement non-exceptional}
Let $X/Z\ni z$ be an $\Rr$-complementary surface that is not exceptional. Then $X/Z\ni z$ has an $n$-complement for some $n\in\{1,2,3,4,6\}$.
\end{lem}
\begin{proof}
There exists an lc but not klt pair $(X/Z\ni z,B)$ such that $K_X+B\sim_{\mathbb R,Z}0$ over a neighborhood of $z$ Let $f: Y\rightarrow X$ be a dlt modification of $(X/Z\ni z,B)$ and let $K_Y+B_Y:=f^*(K_X+B)$. Then $\lfloor B_Y\rfloor\not=0$. By \cite[2.3 Inductive Theorem]{Sho00}, $(Y/Z\ni z,\lfloor B_Y\rfloor)$ has a monotonic $n$-complement $(Y/Z\ni z,B_Y^+)$ for some $n\in\{1,2,3,4,6\}$ for any $z\in Z$, hence $X/Z\ni z$ has an $n$-complement $(X/Z\ni z,f_*B_Y^+)$ for some $n\in\{1,2,3,4,6\}$.
\end{proof}
\begin{cor}\label{cor: surface precise complement relative}
Let $X/Z\ni z$ be an $\Rr$-complementary surface and $\dim Z>0$. Then $X/Z\ni z$ has an $n$-complement for some $n\in\{1,2,3,4,6\}$.
\end{cor}
\begin{proof}
For any pair $(X/Z\ni z,B)$ such that $K_X+B\sim_{\mathbb R,Z}0$ over a neighborhood of $z$, $(X/Z\ni z,B+tf^*z)$ is an lc but not klt pair such that $K_X+B+tf^*z\sim_{\mathbb R,Z}0$ over a neighborhood of $z$, where $t:=\lct(X,B;f^*z)$. The corollary follows from Lemma \ref{lem: complement non-exceptional}.
\end{proof}

\section{The exceptional case}

\subsection{Proof of Theorem \ref{thm: singularity of exceptional fano varieties} and Corollary \ref{cor: 1/42 exceptional fano type surface}}

\begin{lem}\label{lem: bir19 2.50 surface}
Let $d$ be a positive integer, $\Ii\subset [0,1]$ a DCC set, and $\epsilon:=\epsilon(d,\Ii)$. Let $(X/Z,T+aS)$ be a pair such that $X$ is of Fano type over $Z$, $-(K_X+T+aS)$ is nef$/Z$, $T\in \Ii$, $S\not=0$ is a reduced divisor, and $a\in (1-\epsilon,1)$. Then we may run a $-(K_X+T+S)$-MMP$/Z$ which consists of a sequence of divisorial contractions and flips
$$(X,T+S):=(X_0,T_0+S_0)\dashrightarrow (X_1,T_1+S_1)\dashrightarrow\dots\dashrightarrow (X_n,T_n+S_n),$$
such that
\begin{enumerate}
    \item $(X_i,T_i+S_i)$ is lc for each $i$,
    \item $S_n\not=0$, and
    \item $-(K_{X_n}+T_n+S_n)$ is nef$/Z$.
\end{enumerate}
Here $T_i$ and $S_i$ are the strict transforms of $T$ and $S$ on $X_i$ respectively.
\end{lem}
\begin{proof} By Theorem \ref{thm: fano type mds}, we may run a $-(K_X+T+S)$-MMP$/Z$.

(1) Since $-(K_X+T+aS)$ is nef$/Z$, $(X/Z,T+aS)$ is $\mathbb R$-complementary, hence $(X_i/Z,T_i+aS_i)$ is $\mathbb R$-complementary for each $i$. In particular, $(X_i,T_i+aS_i)$ is lc for each $i$. By the definition of $\epsilon$, $(X_i,T_i+S_i)$ is lc for each $i$.

(2) Since $-(K_X+T+aS)$ is nef$/Z$, by the negativity lemma, $X\dashrightarrow X_n$ is $-(K_X+T+aS)$-non-negapositive. Since $X\dashrightarrow X_n$ is a  $-(K_X+T+S)$-MMP, $X\dashrightarrow X_n$ is $-(K_X+T+S)$-negative, hence $X\dashrightarrow X_n$ is $S$-positive, and we get (2).

(3) Suppose not, then this MMP terminates with a $-(K_{X_n}+T_n+S_n)$-Mori fiber space $X_n\rightarrow V$. Then  $-(K_{X_n}+T_n+S_n)$ is anti-ample$/V$. Since $-(K_{X}+T+aS)$ is nef$/Z$, $-(K_{X_n}+T_n+aS_n)$ is nef$/V$, and there exists a real number $c\in [a,1)\subset (1-\epsilon,1)$ such that $K_{X_n}+T_n+cS_n\equiv_V0$. By (1), $(X_n,T_n+cS_n)$ is lc. Let $F$ be a general fiber of $X_n\rightarrow V$, $T_F:=T|_F$, and $S_F:=S|_F$, then $(F,T_F+cS_F)$ is lc and $K_F+T_F+cS_F\equiv 0$. Thus $\epsilon\leq\epsilon_{1}(d,\Ii)\leq\epsilon_1(\dim F,\Ii)\leq 1-c$, a contradiction.
\end{proof} 

\begin{proof}[Proof of Theorem \ref{thm: singularity of exceptional fano varieties}]
Let $a:=\tmld(X,B+G)$ be the total minimal log discrepancy of $(X,B+G)$, and $E$ a divisor over $X$ such that $a(E,X,B+G)=a$. Suppose that $a<\epsilon:=\epsilon(d,\Ii)$. We let $c:=\mult_EG$ and let $e:=1-a-c$. If $E$ is exceptional over $X$, then we let $f: Y\rightarrow X$ be a divisorial contraction which extracts $E$. If $E$ is not exceptional over $X$, then we let $f: Y\rightarrow X$ be the identity morphism. Then we have $K_Y+eE+B_Y=f^*(K_X+B)$, where $B_Y$ is the strict transform of $B$ on $Y$. We let $G_Y:=f^*G$. 

\begin{claim}\label{claim: y is fano type}
$Y$ is of Fano type.
\end{claim}
\begin{proof}
If $Y=X$ then it is clear $Y$ is of Fano type. Otherwise, there exists a klt pair $(X,\Delta)$ such that $-(K_X+\Delta)$ is big and nef. Let $\Delta_Y:=f^{-1}_*\Delta$ and let $a':=a(E,X,\Delta)$, then $a'\leq a<\epsilon<1$, so $(Y,\Delta_Y+(1-a')E)$ is a klt pair such that $-(K_Y+\Delta_Y+(1-a')E)$ is big and nef. Thus $Y$ is of Fano type.
 \end{proof}

\noindent\textit{Proof of Theorem \ref{thm: singularity of exceptional fano varieties} continued}. Since $-(K_Y+(1-a)E+B_Y)\sim_{\mathbb R}G_Y-cE\geq 0$, by Claim \ref{claim: y is fano type} and Theorem \ref{thm: fano type mds}, we may run a $-(K_Y+(1-a)E+B_Y)$-MMP which terminates with a model $T$ such that $-(K_T+(1-a)E_T+B_T)$ is nef, where $E_T,B_T$ are the strict transforms of $E,B$ on $T$ respectively. Since $E$ is not a component of $G_Y-cE$ and the MMP only contracts divisors that are contained in $\Supp(G_Y-cE)$, $E_T\not=0$. By Lemma \ref{lem: bir19 2.50 surface}, we may run a $-(K_T+E_T+B_T)$-MMP which terminates with a model $V$, such that $(V,E_V+B_V)$ is lc, $E_V\not=0$, and $-(K_V+E_V+B_V)$ is nef, where $E_V,B_V$ are the strict transforms of $E_T,B_T$ on $V$ respectively. Since $X$ is of Fano type, $V$ is of Fano type. It is clear that $(V,E_V+B_V)$ is not exceptional.

For any prime divisor $D$ over $X$, we have\begin{align*}
    a(D,X,B)&\geq a(D,Y,(1-a)E+B_Y)\geq a(D,T,(1-a)E_T+B_T)\\
    &\geq a(D,T,E_T+B_T)\geq a(D,V,E_V+B_V).
\end{align*}
By \cite[Lemma 2.17]{Bir19}, $X$ is not exceptional, a contradiction.
\end{proof}

\begin{rem}
The proof of Theorem \ref{thm: singularity of exceptional fano varieties} also works for generalized pairs \cite{BZ16}. For simplicity, we omit the proof.
\end{rem}

\begin{cor}\label{cor: exceptional surface 1/42 lc}
Let $(X,B+G)$ be a pair such that $(X,B)$ is exceptional, $X$ is a Fano type surface, $B\in\{1-\frac{1}{n}\mid n\in\mathbb N^+\}\cup\{1\}$, and $0\leq G\sim_{\mathbb R}-(K_X+B)$. Then $(X,B+G)$ is $\frac{1}{42}$-lc. In particular, $(X,B)$ and $X$ are $\frac{1}{42}$-lc.
\end{cor}
\begin{proof}
It immediately follows from Theorem \ref{thm: singularity of exceptional fano varieties} and Lemma \ref{lem: surface lct 1-gap}.
\end{proof}

\begin{proof}[Proof of Corollary \ref{cor: 1/42 exceptional fano type surface}] It follows from Corollary \ref{cor: exceptional surface 1/42 lc}.
\end{proof}

\begin{cor}\label{cor: index fano exceptional surface}
For any exceptional Fano surface $X$, there exists $I\leq 2\cdot 84^{128\cdot 42^5}$ such that $IK_X$ is Cartier. In particular, $K_X^2\geq\frac{1}{I}$.
\end{cor}
\begin{proof}
It follows from Corollary \ref{cor: 1/42 exceptional fano type surface} and Lemma \ref{lem: index epsilonlc fano surface}.
\end{proof}

\subsection{Exceptional surface complements}

\begin{lem}\label{lem: surface complement contraction}
Let $(X:=\mathbb P^1,B)$ be a pair such that $\deg(K_X+B)\leq 0$ and $B\in\{\frac{k}{12}\mid k\in\mathbb N^+,0\leq k\leq 12\}\cup\{1-\frac{1}{n}\mid n\in\mathbb N^+\}$. Then $(X,B)$ has a monotonic $n$-complement such that $12\mid n$ and $n\leq 276$.
\end{lem}
\begin{proof}
We may write $B=C+D$ where $C,D\geq 0$, $C\in \{1-\frac{1}{n}\mid n\in\mathbb N^+, 12\nmid n\}$, $12D$ is integral, and $C\wedge D=0$. Then the coefficients of $C$ are $\geq\frac{4}{5}$. In particular, $C$ has at most $2$ irreducible components. Possibly adding divisors of the form $\frac{1}{12}p$ to $D$ where $p$ are general points on $X$, we may assume that $0\geq\deg(K_X+B)>-\frac{1}{12}$. We have the following cases.

\medskip

\noindent\textbf{Case 1}. $C=0$. Then $(X,B)$ is a $12$-complement of itself.

\medskip

\noindent\textbf{Case 2}.  $C$ has $1$ irreducible component $C_1$. Then $C=aC_1$ for some $a\in (0,1)$ and $(X,\frac{\lceil 12a\rceil}{12}C_1+D)$ is a monotonic $12$-complement of $(X,B)$.

\medskip

\noindent\textbf{Case 3}. $C$ has $2$ irreducible components $C_1,C_2$. We have $C=a_1C_1+a_2C_2$. Possibly switching $C_1,C_2$, we may assume that $a_1\leq a_2$. If  $D=0$, then $(X,C_1+C_2)$ is a monotonic $1$-complement of $(X,B)$. If $D\not=0$, then $a_1\leq\frac{23}{24}$. Let $m$ be the denominator of $a_1$, then $m\leq 24$ and $(X,a_1C_1+(2-\deg D-a_1)C_2+D)$ is a monotonic $\lcm(12,m)$-complement of $(X,B)$. Since $12\mid\lcm(12,m)$ and $\lcm(12,m)\leq 276$, we are done.
\end{proof}

\begin{thm}\label{thm: surface complement exceptional}
Let $X$ be an $\Rr$-complementary exceptional surface.
\begin{enumerate}
    \item If $\kappa(-K_X)=0$, then $X$ has an $n$-complement for some $n\leq 21$.
    \item If $\kappa(-K_X)=1$, then $X$ has an $n$-complement for some $n$ such that $12\mid n$ and $n\leq 276$.
    \item If $\kappa(-K_X)=2$, then $X$ has an $n$-complement for some $n\leq 192\cdot 84^{128\cdot 42^5}$.
\end{enumerate}
\end{thm}
\begin{proof}
There exists a klt pair $(X,B)$ such that $K_X+B\sim_{\mathbb R}0$. Thus  $(X,(1+\delta)B)$ is klt for some $0<\delta\ll 1$, so we may run a $(K_X+(1+\delta)B)$-MMP which terminates with good minimal model $X'$. Since $K_X+(1+\delta)B\sim_{\mathbb R}-\delta K_X$, this is also a $-K_X$-MMP. By abundance for klt surfaces, $-K_{X'}$ is semi-ample. Possibly replacing $X$ with $X'$, we may assume that $-K_{X}$ is semi-ample. 

If  $\kappa(-K_X)=0$, then $K_X\equiv 0$. Hence $nK_X\sim 0$ for some positive integer $n\leq 21$ \cite{Bla95,Zha91,Zha93}, and $X$ is an $n$-complement of itself for some $n\leq 21$.

If $\kappa(-K_X)=1$, then $-K_X$ defines a contraction $f: X\rightarrow Z$. By Kodaira's canonical bundle formula, we have
$$12K_X\sim 12f^*(K_Z+B_Z+M_Z)$$
such that $12M_Z$ is an integral divisor, $B_Z\in\{1-\frac{1}{n}\mid n\in\mathbb N^+\}$, and $$B_Z=\sum_{z\in Z}(1-\lct(X,0;f^*z))z.$$
We may choose $M_Z$ such that $B_Z\wedge M_Z=0$ and $(Z,B_Z+M_Z)$ is klt. Since $-K_X$ is semi-ample, $\deg(K_Z+B_Z+M_Z)\leq 0$, hence $Z$ is either an elliptic curve or $\mathbb P^1$. If $Z$ is an elliptic curve, then $B_Z=M_Z=0$ and $12K_X$ is base-point-free, hence $X$ is a $12$-complement of itself. If $Z=\mathbb P^1$, by Lemma \ref{lem: surface complement contraction}, there exists an integer $n\leq 276$ such that $12\mid n$ and $(Z,B_Z+M_Z)$ has a monotonic $n$-complement $(Z,B_Z+G+M_Z)$. By the construction of $B_Z$, $(X,f^*G)$ lc, hence $(X,f^*G)$ is an $n$-complement of $X$.

If $\kappa(-K_X)=2$, then $-K_X$ defines a birational morphism $f: X\rightarrow Y$. We have $K_X=f^*K_Y$. Possibly replacing $X$ with $Y$, we may assume that $X$ is Fano. By Corollary \ref{cor: 1/42 exceptional fano type surface}, $X$ is $\frac{1}{42}$-lc. By Lemma \ref{lem: index epsilonlc fano surface}, $IK_X$ is Cartier for some $I\leq 2\cdot 84^{128\cdot 42^5}$. By the effective base-point-freeness theorem (\cite[Theorem 1.1, Remark 1.2]{Fuj09} \cite[1.1 Theorem]{Kol93}), $|-96IK_X|$ is base-point-free. In particular, $X$ has a $96I$-complement. 
\end{proof}

\section{Proof of the main theorems}

\begin{proof}[Proof of Theorem \ref{thm: surface precise complement}]
If $\dim Z>0$ or $\dim Z=0$ and $X$ is not exceptional, the theorem follows from Lemma \ref{lem: complement non-exceptional} and Corollary \ref{cor: surface precise complement relative}. Otherwise, $\dim Z=0$ and $X$ is exceptional, and the theorem follows from Theorem \ref{thm: surface complement exceptional}.
\end{proof}

\begin{proof}[Proof of Corollary \ref{cor: upper bound alpha invariant surface}]
Recall that Tian's $\alpha$-invariant for a variety $X$ is defined as
$$\alpha(X):=\inf\{t\geq 0\mid \lct(X,0;D)\mid D\in |-K_X|_{\mathbb Q}\}.$$
We may assume that $\alpha(X)>1$. If $\kappa(-K_X)\leq 1$, then by Theorem \ref{thm: surface complement exceptional},  $X$ has an $n$-complement $(X,G)$ for some $n\leq 276$, hence $(X,nG)$ is not klt, and $\alpha(X)\leq 276$. Thus we may assume that $\kappa(-K_X)=2$. We may run a $(-K_X)$-MMP and replace $X$ with the canonical model of $-K_X$, and assume that $-K_X$ is ample. By Corollary \ref{cor: index fano exceptional surface}, $K_X^2\geq\frac{1}{I}$. Since $\alpha(X)^2\cdot\vol(-K_X)\leq 9$ (cf. \cite[Theorems A,D]{BJ20}), $\alpha(X)\leq 3\sqrt{I}\leq 3\sqrt{2}\cdot 84^{64\cdot 42^5}$.
\end{proof}

\section{Further remarks}

\begin{rem}[Reasonable and optimal bounds]
Recall that we expect $\epsilon_1(2,\{0\})=\frac{1}{13}$. If we can prove this, then the bound of $n$ in Theorem \ref{thm: surface precise complement} can be improved to $192\cdot 26^{128\cdot 13^5}\approx 10^{10^{7.8}}$. This is much smaller than the current bound, albeit it is still expected to be far from optimal. On the other hand, if one can get a better bound for $I=I(\epsilon)$ in Lemma \ref{lem: index epsilonlc fano surface}, then the bound of $n$ may be greatly improved. For example, \cite[Proof of Lemma 4.9]{LX22} actually implies that the local Cartier index of any $\frac{2}{5}$-klt weak Fano surface is $\leq 19$. With a little more effort, one can show that the global Cartier index of any $\frac{2}{5}$-klt weak Fano surface is $\leq 385$. This is much smaller than the bound given by Lemma \ref{lem: index epsilonlc fano surface} which is $2\cdot 5^{12500}$. By applying the arguments in this paper, we shall get $n\leq 36960$ and $\alpha(X)\leq 3\sqrt{385}\approx 58.86$ for exceptional $\frac{2}{5}$-klt surfaces.
\end{rem}

\begin{rem}[Explicit bound for pairs]\label{rem: explicit bound for pairs}
One may also ask whether we can find an explicit bounded of $n$ for $n$-complements of surface pairs $(X,B)$. 

For pairs with finite rational coefficients, the bound is computable via the methods introduced in \cite{AM04}, but may be much larger than the case when $B=0$. This is because the bound $\epsilon(d,\Ii)=\frac{1}{42}$ in Lemma \ref{lem: surface lct 1-gap} will be changed to a number which is very close to $1$ as in \cite[3.5]{AM04} when $\Ii$ is not the standard set. It is very difficult to represent that number in a very explicit function of the common denominator of the coefficient set (even when the coefficient set is $\{0,\frac{1}{3}\}$, for example). We also need to go through the inductive arguments as in \cite[2.3 Inductive Theorem]{Sho00} for non-exceptional complements. 

For pairs with finite (maybe irrational) coefficients or DCC coefficients, a theory on ``explicit uniform rational polytopes" (cf. \cite{HLS19}) is needed, which is still unknown. 

For pairs with coefficients in $[0,1]$, one needs to go through all the previous simpler cases and check the details of the proof of \cite[Theorem 3]{Sho20} and avoid using any inexplicit boundedness result. This is considered to be much more difficult. See \cite{CHX23} for a similar result.
\end{rem}

\begin{rem}[Explicit bound of threefold mlds]
By Remark \ref{rem: explicit bound for pairs} and following the details of the proof of \cite{HLL22}, one will be able to provide a computable lower bound of the $1$-gap of threefold mlds for pairs with finite rational coefficients (or, more generally, hyperstandard rational coefficients). This is because all other constants in the proof of \cite{HLL22} can be explicitly bounded except the $t$ in \cite[Lemma 6.4]{HLL22}. Here an explicit boundedness of $n$-complement for threefold singularities is needed, but this just follows from Theorem \ref{thm: surface precise complement} and Remark \ref{rem: explicit bound for pairs}. Nevertheless, such bound will, again, be far from being optimal. For example, when we have standard coefficients, the $1$-gap is expected to be $\frac{1}{42}$ (by Remark \ref{rem: gloabl local correspondenc mld gap} and Lemma \ref{lem: surface lct 1-gap}), but we can only show that the $1$-gap is $\geq\epsilon$ for some $\epsilon\approx\frac{1}{192\cdot 42^{128\cdot 42^5}}$.
\end{rem}

\end{document}